# Generating functions for the generalized Li's sums

**S. K. Sekatskii** (LPMV, Ecole Polytechnique Fédérale de Lausanne, Switzerland)

Recently, in arXiv:1304.7895; *Ukrainian Math. J.* – 2014.- **66**. - P. 371 – 383, we presented the generalized Li's criterion. This is the statement that the sums $\lambda_{n,b,\sigma} = \sum_{\rho} (1 - \left(\frac{\rho+b}{\rho-b-2\sigma}\right)^n)$, taken over all Riemann xi-function zeroes taking into account their multiplicity (complex conjugate zeroes are to be paired when summing whenever necessary) for any $n=1, 2, 3...$ and any real $b > -\sigma$, are non-negative if and only if there are no Riemann function zeroes with $\text{Re}\,\rho > \sigma$. For any $n=1, 2, 3...$ and any real $b < -\sigma$, these sums are non-negative if and only if there are no Riemann function zeroes with $\text{Re}\,\rho < \sigma$; correspondingly, for $\sigma = 1/2$ and $b \neq -1/2$ such non-negativity is equivalent to the Riemann hypothesis. In this Note we obtain generation functions for this generalized criterion demonstrating the Taylor expansion ($b \neq -\sigma$):

$$\ln(\xi\left(b + 2\sigma + \frac{(2b+2\sigma)z}{1-z}\right)) = \ln(\xi(b+2\sigma)) + \sum_{n=1}^{\infty} \frac{\lambda_{n,b,\sigma}}{n} z^n .$$

**Introduction**

In paper [1], Li has established the following criterion (now bearing his name) equivalent to the Riemann hypothesis (RH):

**Theorem 1 (Li's criterion).** *Riemann hypothesis is equivalent to the non-negativity of the following sums over the Riemann zeta-function non-trivial zeroes $\rho$ taking into account their multiplicity (complex conjugate zeroes are to be paired whenever necessary) for all positive integers n:*

$$\lambda_n := \sum_{\rho} (1 - (\frac{\rho}{\rho - 1})^n) \geq 0 \qquad (1).$$

(See e.g. [2] for standard definitions and discussion of the general properties of the Riemann zeta-function). In the same paper, Li, by establishing the Taylor expansion

$$\ln(\xi\left(\frac{z}{z-1}\right)) = -\ln 2 + \sum_{n=1}^{\infty} \frac{\lambda_n}{n} z^n \qquad (2),$$

also found the generating function for his sums. (Li gave Taylor expansion in terms of the logarithmic derivative $\frac{\xi'}{\xi}\left(\frac{z}{z-1}\right) = \sum_{n=0}^{\infty} \lambda^{n+1} z^n$; actually, expansion at question has been established earlier by Keiper [3]).

Recently, we have generalized Li's criterion [4], and for our current purposes we combine a few previous results into the following theorem.

**Theorem 2.** *For any real $\sigma$ and any real $b > -\sigma$, the following sums over the Riemann zeta-function non-trivial zeroes $\rho$ taking into account their multiplicity (complex conjugate zeroes are to be paired whenever necessary) are non-negative for all positive integers n,*

$$\lambda_{n,b,\sigma} := \sum_{\rho} (1 - (\frac{\rho + b}{\rho - b - 2\sigma})^n) \geq 0 \qquad (3),$$

*if and only if there are no non-trivial Riemann zeta-function zeroes with* $\operatorname{Re}\rho > \sigma$. *For any real* $b < -\sigma$, *this takes place if and only if there are no non-trivial Riemann zeta-function zeroes with* $\operatorname{Re}\rho < \sigma$. *Correspondingly, the non-negativity of the sums* $\lambda_{n,b,1/2} := \sum_{\rho}(1-(\frac{\rho+b}{\rho-b-1})^n) \geq 0$ *for all n and all* $b \neq -1/2$ *is equivalent to the Riemann hypothesis.*

In relation with these theorems, one can immediately put forward a question concerning the construction of generating functions for the generalized sums $\lambda_{n,b,\sigma}$. An answer to this question is the aim of the present short Note.

**2. Generating functions for the sums** $\lambda_{n,b,\sigma}$

**Theorem 3.** *For any real* $b$, $\sigma$, $\sigma \neq -b$, *the sums* $\lambda_{n,b,\sigma} := \sum_{\rho}(1-(1+\frac{2b+2\sigma}{\rho-b-2\sigma})^n) = \sum_{\rho}(1-(\frac{\rho+b}{\rho-b-2\sigma})^n)$, *which appear in the formulation of the generalized Li's criterion, are given by the coefficients of the following Taylor expansion:*

$$\ln(\xi\left(b+2\sigma+\frac{(2b+2\sigma)z}{1-z}\right)) = \ln(\xi(b+2\sigma)) + \sum_{n=1}^{\infty}\frac{\lambda_{n,b,\sigma}}{n}z^n \qquad (4).$$

**Proof.** Proof is a straightforward generalization of approach given in [1]. We start from the well known property concerning the representation of the Riemann *xi*-function as an infinite Hadamard product over its zeroes [2]

$$\xi(s) = \frac{1}{2}\prod_{\rho}(1-\frac{s}{\rho}) \qquad (5)$$

and then introduce the function $\xi^*(z) := \xi(z+b+2\sigma)$. We have

$$\xi^*(s) = \xi(b+2\sigma) \cdot \prod_{\rho^*}(1-\frac{s}{\rho^*}) \tag{5a}$$

where $\rho^* = \rho - b - 2\sigma$ are zeroes of $\xi^*(z)$. (This, similarly to (5), is a direct consequence of the general Hadamard product formula and the relation $\frac{\xi'}{\xi}(b+2\sigma) = \sum_{\rho}\frac{1}{b+2\sigma-\rho}$ [2]; complex conjugate zeroes are to be paired when summing). Now we take $s = -\frac{(2b+2\sigma)z}{z-1}$ whence $1-\frac{s}{\rho^*} = \frac{1-z(1+\frac{2b+2\sigma}{\rho^*})}{1-z}$.

We have, by taking a logarithm and applying Taylor expansion

$$\ln(1-\frac{s}{\rho^*}) = \ln(1-z(1+\frac{2\sigma+2b}{\rho^*})) - \ln(1-z) = \sum_{n=1}^{\infty}\frac{1}{n}(1-(1+\frac{2b+2\sigma}{\rho^*})^n)z^n$$

so that, according to eq. (5a), the sums $\lambda_{n,b,\sigma} := \sum_{\rho}(1-(1+\frac{2b+2\sigma}{\rho-b-2\sigma})^n)$ at question are given by Taylor expansion of the function

$$\ln(\xi^*(\frac{(2b+2\sigma)z}{1-z})) = \ln(\xi(b+2\sigma+\frac{(2b+2\sigma)z}{1-z})) = \ln(\xi(\frac{b+2\sigma+bz}{1-z}))$$ at the point $z=0$:

$$\ln(\xi\left(b+2\sigma+\frac{(2b+2\sigma)z}{1-z}\right)) = \ln(\xi(b+2\sigma)) + \sum_{n=1}^{\infty}\frac{\lambda_{n,b,\sigma}}{n}z^n.$$

## 3. Concluding remarks

In particular, if we take $b = -2\sigma$ we have the function $\ln(\xi(\frac{2\sigma z}{z-1}))$ which is quite similar to the original Li's function: the non-negativity of all coefficients of the expansion $\ln(\xi\left(\frac{2\sigma z}{z-1}\right)) = -\ln 2 + \sum_{n=1}^{\infty}\frac{\lambda_{n,-2\sigma,\sigma}}{n}z^n$ for some $0 < \sigma < 1/2$ is equivalent to the statement that there are no zeroes with $\text{Re}\,\rho < \sigma$. Is it possible that an additional factor $\sigma$ occurring here might be a useful leverage for an analysis of this expansion by comparing it with some other relevant expansions, e. g. $\ln(\xi(z))$, similarly to works of Coffey [4 - 6] and others [7]?